\documentclass[onecolumn,hideversion]{cdcarticle}

\usepackage{amsmath}
\usepackage{amssymb}
\usepackage{amsthm}
\usepackage{graphicx}

\usepackage{stdmath}

\graphicspath{{graphics_included/}}

\def\bk{\hspace{-0.75mm}}

\def\S{\mathcal{S}}
\def\X{\mathcal{X}}
\def\U{\mathcal{U}}

\newcommand{\argmin}[1]{\hbox{$\underset{#1}{\mbox{argmin}}$}}
\newcommand{\nchoosek}[2]{\left(\begin{array}{c}#1\\#2\end{array}\right)}

\begin{document}

\title{On Decentralized Policies for the \\ Stochastic $k$-Server Problem}

\author{%
  Randy Cogill\footnotesymbol{1}
  \and
  Sanjay Lall\footnotesymbol{2}
  }

\note{}

\maketitle

%%%%%%%%%%%%%%%%%%%%%%%%%%%%%%%%%%%%%%%%%%%%%%%%%%%%%%%%%%%%%%%%%%%%%%%%%%%%%%
% the footnotes

\makefootnote{1}{Department of Electrical Engineering,
  Stanford University, \\
  Stanford, CA 94305, U.S.A.
  Email: rcogill@stanford.edu}

\makefootnote{2}{
  Department of Aeronautics and Astronautics,
  Stanford University, Stanford CA 94305-4035, U.S.A.
  Email lall@stanford.edu}

\makefootnote{1}{The first author was partially supported
  by a Stanford Graduate Fellowship.}

\makefootnote{1,2}{Partially supported by the Stanford URI
  \emph{Architectures for Secure and Robust Distributed
    Infrastructures}, AFOSR DoD award number 49620-01-1-0365.  }

%%%%%%%%%%%%%%%%%%%%%%%%%%%%%%%%%%%%%%%%%%%%%%%%%%%%%%%%%%%%%%%%%%%%%%%%%%%%%%
\begin{abstract}

In this paper we study a dynamic resource allocation problem which we call the \emph{stochastic $k$-server problem}. In this problem, requests for some service to be performed appear at various locations over time, and we have a collection of $k$ mobile servers which are capable of servicing these requests. When servicing a request, we incur a cost equal to the distance traveled by the dispatched server. The goal is to find a strategy for choosing which server to dispatch to each incoming request which keeps the average service cost as small as possible.

In the model considered in this paper, the locations of service requests are drawn according to an IID random process. We show that, given a statistical description of this process, we can compute a simple decentralized state-feedback policy which achieves an average cost within a factor of two of the cost achieved by an optimal state-feedback policy. In addition, we demonstrate similar results for several extensions of the basic stochastic $k$-server problem.

\end{abstract}

%%%%%%%%%%%%%%%%%%%%%%%%%%%%%%%%%%%%%%%%%%%%%%%%%%%%%%%%%%%%%%%%%%%%%%%%%%%%%%

\section{Introduction}

Recently, there has been great interest in the study of coordination strategies for teams of Unmanned Aerial Vehicles (UAVs). In particular, many researchers have focused on methods for designing efficient mission plans, under which a series of tasks can be carried out by a team of vehicles. A common high-level formulation of this type of problem consists of a series of waypoints that must be visited by the vehicles, with the goal of designing a strategy for visiting each of the waypoints in a manner which minimizes some measure of the overall travel time. When the set of locations to visit is known ahead of time, it is possible to plan the mission offline, and each vehicle can perform its own tasks without requiring communication among the vehicles \cite{Gil03,Richards02}. In a dynamic environment, where waypoints may appear as the system is operation, the mission cannot be planned entirely ahead of time. Such a formulation is considered in \cite{Frazzoli04}. However, due to limited computational and communication resources, it is generally not feasible to consider coordination strategies which require complete communication among the vehicles during system operation. The general problem considered in this paper is motivated by the problem of multi-vehicle coordination in a dynamic environment.

The well known $k$-server problem is a natural model for dynamic task assignment problems with distance-based costs. Roughly speaking, the $k$-server problem is as follows. We are given a set of locations, and requests for services to be performed originate sequentially from these locations. We have a collection of $k$ mobile servers which are capable of servicing these requests. At each point in time we must choose a server to serve the current request, and we incur a cost equal to the distance traveled by the dispatched server. The goal is to find a strategy for choosing which servers to dispatch to each incoming request which keeps the average service cost as small as possible.

The $k$-server problem has been well studied for the problem formulation where the demand sequence may be arbitrary. Most of the literature on the $k$-server problem has focused on the competitive analysis of online algorithms. An online algorithm is a strategy which makes decisions based only on the knowledge of present and past requests, and competitive analysis seeks to compare the performance of specific online algorithms with the performance of an optimal strategy which knows the entire request sequence. The best known results for the $k$-server problem show that a particular online algorithm (which requires intensive computation to implement) achieves an overall cost which is essentially within a factor of $2k-1$ of optimal \cite{Koutsoupias94}. The reader is referred to \cite{Floratos97} for a survey on the $k$-server problem, online algorithms and competitive analysis.

In this paper we consider a variation of the $k$-server problem, where the locations of the service requests are drawn at random according to an IID random process. With this stronger assumption on the demand sequence, it is possible to show that a simple, practical strategy can achieve performance comparable to an optimal state feedback strategy. Specifically, we show that, given a statistical description of the request sequence, we can compute a simple decentralized state-feedback policy achieves an overall cost within a factor of two of the cost achievable by an optimal state-feedback policy. A decentralized policy has the property that, once the policy is determined, no communication between the servers is required for its implementation. In addition, we demonstrate similar results for several extensions of the basic stochastic $k$-server problem.

%%%%%%%%%%%%%%%%%%%%%%%%%%%%%%%%%%%%%%%%%%%%%%%%%%%%%%%%%%%%%%%%%%%%%%%%%%%%%%

\section{Problem formulation}

In this section we give a precise formulation of the stochastic $k$-server problem. In our formulation, the servers are positioned and the requests originate at points in some finite set $\S$. The set $\S$ is equipped with a metric $d:\S\times \S\rightarrow \mathbb{R}_+$. At each time step $t\in \mathbb{Z}_+$, service at some point $\widehat{x}(t) \in \S$ is requested and the $k$ servers reside at the points $x_1(t),\ldots,x_k(t)\in\S$. Exactly one server must be chosen to service the request at $\widehat{x}(t)$. If server $u(t)$ is chosen, then a cost of $d(x_{u(t)}(t),\widehat{x}(t))$ is incurred and server $u(t)$ is relocated to the point $\widehat{x}(t)$. That is, $x_{u(t)}(t+1) = \widehat{x}(t)$ and $x_i(t+1) = x_i(t)$ for all $1 \le i \le k$ such that $i\ne u(t)$. The next service request $\widehat{x}(t+1)$ is then randomly chosen. In our model, each $\widehat{x}(t)$ is drawn according to the probability mass function $p:\S\rightarrow[0,1]$, and is independent of $\widehat{x}(\tau)$ for all $\tau\ne t$. The goal of the problem is to determine a strategy for assigning servers to service requests which keeps the average cost incurred in each time step as small as possible. This is illustrated in Figure 1.

\vspace{2mm}

\begin{figure}[ht!]
\label{fig:k_server_1}
\begin{center}
\scalebox{0.5}{\includegraphics{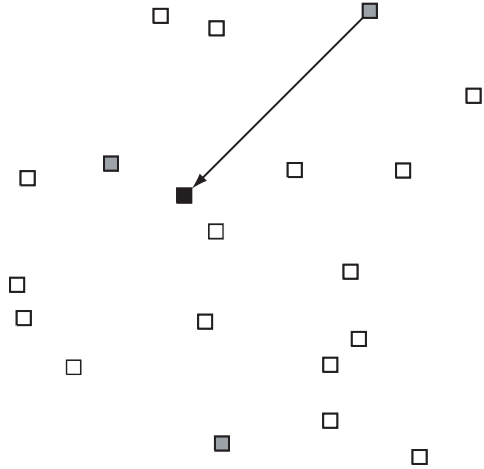}}
\caption{Gray squares represent server locations and the black square represents the location of the service request. One server will move to the request location and incur a cost equal to the distance traveled.}
\end{center}
\end{figure}

The problem described in the previous paragraph can be formulated as a finite state Markov decision process with average cost criteria (see, for example, \cite{Puterman94}). In general, finite state Markov decision processes have a finite state space $\X$, and a finite set of actions $\U$ available at each time step. Taking action $u \in\U$ when in state $y \in\X$ incurs a cost $r(y,u)$. After taking action $u$ in state $y$, the system state in the next time period is $x \in\X$ with probability $\textbf{Pr}(X(t+1)=x~|~X(t)=y,U(t)=u)$.

A static state-feedback policy is a decision rule in which each $u(t)$ is chosen according to a function $\mu:\X\rightarrow \U$ of the current state $x(t)$. The steady-state average per-period cost under the policy $\mu$ is
\begin{eqnarray*}
\lefteqn{J(\mu,x(0)) =} \\ && \lim_{t\rightarrow\infty}\frac{1}{t+1}\sum_{k=0}^t E[r(X(k),\mu(X(k)))~|~X(0)=x(0)].
\end{eqnarray*}
We denote a policy which minimizes this cost by $\mu_*$.

The obvious formulation of the stochastic $k$-server problem as a Markov decision process has the state at time $t$ given by $x(t) = (\widehat{x}(t),x_1(t),\ldots,x_k(t))$, the current service request location together with the set of current server locations. The state space is as a subset of $\S^{k+1}$ since we may exclude, without loss of generality, all states which have more than one server assigned to a particular location. The action $u(t)$ taken at time $t$ is the index of the chosen server, and the action space is $\U = \{1,\ldots,k\}$. The cost incurred at time $t$ is $r(x(t),u(t)) = d(x_{u(t)}(t),\widehat{x}(t))$, the distance from the dispatched server to the current service request. Under a static state-feedback policy, the state evolves according to a Markov chain since $\widehat{x}$ is an IID random process and for each $t$, $x_1(t),\ldots,x_k(t)$ depends only on the previous state.

Although algorithms exist for determining an optimal state-feedback policy for average cost Markov decision processes, they are generally not practical for this problem. One reason is that, under the formulation above, the system has $|\S|^{k+1}$ discrete states. Numerical computation of an optimal policy will be intractable even for relatively small values of $|\S|$ and $k$. Also, even if the optimal policy could be computed, this policy may not lend itself to practical implementation. In particular, the optimal policy may be structured so that the decision $u(t)$ must be made based on the knowledge of all server locations at time $t$. This means that all servers would be required to communicate their current locations to all other servers before each decision could be made. In the next section, we will show that a fairly simple decentralized strategy can achieve an average per-period cost within a factor of two of an optimal centralized strategy.

%%%%%%%%%%%%%%%%%%%%%%%%%%%%%%%%%%%%%%%%%%%%%%%%%%%%%%%%%%%%%%%%%%%%%%%%%%%%%%

\section{Main result}

In this section we will consider decentralized policies for the $k$-server problem. After introducing decentralized policies, we will show that there is a decentralized policy that can achieve performance close to that of an optimal policy.

In a general state feedback policy, the decision of which server to dispatch to a request depends on the location of the request as well as the current location of all servers. In contrast, a decentralized policy is a policy in which each server makes a decision to serve the current request without knowledge of the locations of other servers. Given that one and only one server must respond to each request, it is necessary that decentralized policies have a special `partition' structure. That is, decentralized policies partition the set $\S$ into $k$ disjoint sets $\S_1,\ldots,\S_k$, and server $i$ serves location $\widehat{x}$ if and only if $\widehat{x} \in \S_i$. This is illustrated in Figure 2. 

\vspace{2mm}

\begin{figure}[ht!]
\begin{center}
\scalebox{0.5}{\includegraphics{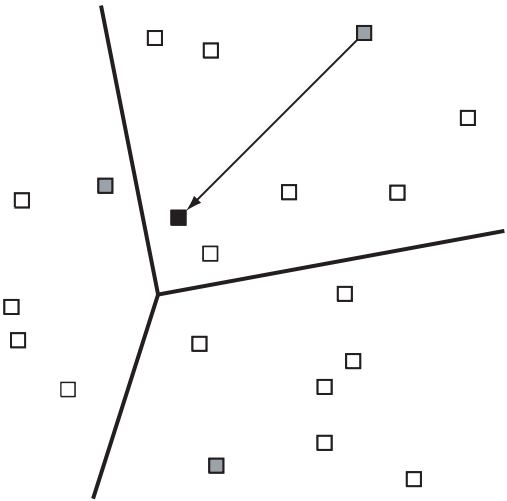}}
\caption{Locations are separated into disjoint partitions. Servers only serve locations in their partition.}
\end{center}
\label{fig:k_server_2}
\end{figure}

\noindent It turns out that there is always a decentralized policy for any instance of the stochastic $k$-server problem that can achieve an average cost comparable to the optimal centralized cost. This policy, which we will call $\mu_d$, is constructed as follows:

\begin{enumerate}
\item Compute the $m^*_1,\ldots,m^*_k$ which minimize
\[
\sum_{s\in \S}p(s)\min_{i\in\U}\{d(m_i,s)\}.
\]
\item Construct the disjoint partitions $\S_1,\ldots,\S_k$, where
\[
\S_i = \{s\in\S ~|~ d(m^*_i,s) \le d(m^*_j,s) \text{ for } j=1,\ldots,k\}.
\]
\item Let $\mu_d(x) = i$ if $\widehat{x} \in \S_i$.
\end{enumerate}

\noindent Performance of this policy relative to an optimal policy is characterized in the following theorem, which is the main result of this paper.

\begin{thm}
\label{thm:factor_two}
The cost of the decentralized policy $\mu_d$ satisfies
\[
J(\mu_d,x(0)) \le 2J(\mu_*,x(0))
\]
for all $x(0) \in \X$.
\end{thm}

\noindent In order to prove Theorem \ref{thm:factor_two}, we will employ a result which allows one to generate performance bounds for general Markov decision processes. This result is proven in \cite{Cogill06} for the case of general measurable state spaces, and is presented here for the finite state space case. 

\vspace{1mm}

\begin{lem}
\label{lem:markov_bounds}
Consider a finite state Markov decision process with average cost criteria. For any state feedback policy $\mu:\X\rightarrow\U$ and any function $h_U:\X\rightarrow\mathbb{R}$,
\[
J(\mu,x(0)) \le \sup_{x\in\X}\{r(x,\mu(x))+\Delta_U(x)\}
\]
for all $x(0)\in \X$, where 
\[
\Delta_U(x)\bk =\bk E[h_U(X(t+1))~|X(t)\bk=\bk x,U(t)\bk=\bk \mu(x)]-h_U(x).
\]
Moreover, for any function $h_L:\X\rightarrow\U$,
\[
J(\mu_*,x(0)) \ge \inf_{x\in\X,u\in\U}\{r(x,u)+\Delta_L(x,u)\}
\]
for all $x(0)\in\X$, where 
\[
\Delta_L(x,u)\bk =\bk E[h_L(X(t+1))~|X(t)\bk=\bk x,U(t)\bk=\bk u]-h_L(x).
\]
\end{lem}

\vspace{1mm}

\noindent Given the result in Lemma~\ref{lem:markov_bounds}, we can now prove Theorem~\ref{thm:factor_two}.

\vspace{1mm}

\noindent \textbf{Proof of Theorem \ref{thm:factor_two}.} First we will find a lower bound on $J(\mu_*,x(0))$ using Lemma \ref{lem:markov_bounds} with
\[
h_L(x) = \min_{i\in\U}\{d(x_i,\widehat{x})\}.
\]
For this choice of $h_L$,
\[
\Delta_L(x,u) = 
\sum_{s\in\S}p(s)\min_{i\in\U}\{d(x_i(t+1),s)\}-\min_{i\in\U}\{d(x_i,\widehat{x})\},
\]
where
\begin{eqnarray*}
x_i(t+1) = \left\{
\begin{array}{cl}
\widehat{x} & \text{if } i = u \\
x_i & \text{otherwise}
\end{array}
\right. .
\end{eqnarray*}
Since $d(x_u,\widehat{x}) \ge h_L(x)$ for all $x\in \X$ and $u\in\U$, by Lemma \ref{lem:markov_bounds} we have
\begin{eqnarray}
\label{eqn:k_median_cost}
J(\mu_*,x(0)) \ge \min_{m\in\S^k}\left\{\sum_{s\in\S}p(s)\min_{i\in\U}\{d(m_i,s)\} \right\}
\end{eqnarray}
for all $x(0) \in \X$.

Let $m^*$ denote the minimizing $m$ in (\ref{eqn:k_median_cost}). Recall that the decentralized policy $\mu_d$ divides the set $\S$ into disjoint partitions $\S_1,\ldots,S_k$, where
\[
\S_i = \{s\in\S ~|~ d(m^*_i,s) \le d(m^*_j,s) \text{ for } j=1,\ldots,k\}.
\]
We will find an upper bound on $J(\mu_d,x(0))$ using Lemma \ref{lem:markov_bounds} with
\[
h_U(x) = 2\min_{i\in\U}\{d(m^*_i,\widehat{x})\}+\sum_{i=1}^kd(x_i,m^*_i).
\]
For this choice of $h_U$,
\begin{eqnarray*}
\lefteqn{r(x,\mu_d(x)) + \Delta_U(x) =} \\ && 2\sum_{s\in\S}p(s)\min_{i\in\U}\{d(m^*_i,s)\} + \\ && d(x_{\mu_d(x)},\widehat{x}) - d(\widehat{x},m^*_{\mu_d(x)}) - d(x_{\mu_d(x)},m^*_{\mu_d(x)}).
\end{eqnarray*}
Since $d$ is a metric,
\[
d(x_{\mu_d(x)},\widehat{x}) \le d(x_{\mu_d(x)},m^*_{\mu_d(x)}) + d(m^*_{\mu_d(x)},\widehat{x}),
\]
and therefore
\begin{eqnarray*}
J(\mu_d,x(0)) &\le& 2\sum_{s\in\S}p(s)\min_{i\in\U}\{d(m^*_i,s)\} \\
&\le& 2J(\mu_*,x(0)).
\end{eqnarray*}
$\blacksquare$

%%%%%%%%%%%%%%%%%%%%%%%%%%%%%%%%%%%%%%%%%%%%%%%%%%%%%%%%%%%%%%%%%%%%%%%%%%%%%%

\section{Computing decentralized policies}

It was shown in the last section that finding a decentralized policy which achieves an average cost within a factor of two of optimal reduces to finding the $m^*_1,\ldots,m^*_k$ minimizing
\[
\sum_{s\in\S}p(s)\min_{i\in\U}\{d(m_i,s)\}.
\]
In other words, a decentralized policy for our dynamic problem can be determined by solving a static combinatorial optimization problem. This static problem has been well studied, and is known as the \emph{$k$-median problem}.

The number of possible solutions to the $k$-median problem is $\nchoosek{|\S|}{k}$. Unfortunately, there are no known algorithms for finding an optimal solution with computational requirements that scale well with $k$. However, much study has been devoted to efficient approximation algorithms for this problem. In this section we will show that the result of the previous section can be combined with known results on approximation algorithms for the $k$-median problem to obtain efficient algorithms for computing decentralized policies for the stochastic $k$-server problem. 

Suppose $\widetilde{m}_1,\dots,\widetilde{m}_k$ is a suboptimal solution to the $k$-median problem. Let $\widetilde{\mu}_d$ be the decentralized policy constructed with the disjoint partitions $\S_1,\ldots,\S_k$, where
\[
\S_i = \{s\in\S ~|~ d(\widetilde{m}_i,s) \le d(\widetilde{m}_j,s) \text{ for } j=1,\ldots,k\}.
\]
The following lemma relates the performance of the policy $\widetilde{\mu}_d$ to the quality of the suboptimal $k$-median solution $\widetilde{m}_1,\dots,\widetilde{m}_k$.

\begin{lem}
\label{lem:subopt_k_median}
Suppose the suboptimal $k$-median solution $\widetilde{m}_1,\dots,\widetilde{m}_k$ satisfies
\[
\sum_{s\in\S}p(s)\min_{i\in\U}\{d(\widetilde{m}_i,s)\} \le
\rho \sum_{s\in\S}p(s)\min_{i\in\U}\{d(m^*_i,s)\}
\]
for some $\rho \ge 1$. Then 
\[
J\bigl(\widetilde{\mu}_d,x(0)\bigr) \le 2\rho J(\mu_*,x(0))
\]
for all $x(0) \in \X$.
\end{lem}

\vspace{1mm}

\begin{proof}
We can find an upper bound on $J\bigl(\widetilde{\mu}_d,x(0)\bigr)$ using Lemma \ref{lem:markov_bounds} with
\[
h_U(x) = 2\min_{i\in\U}\{d(\widetilde{m}_i,\widehat{x})\}+\sum_{i=1}^kd(x_i,\widetilde{m}_i).
\]
Proceeding exactly as in the proof of Theorem \ref{thm:factor_two}, we obtain
\begin{eqnarray*}
J\bigl(\widetilde{\mu}_d,x(0)\bigr) &\le& 2 \sum_{s\in\S}p(s)\min_{i\in\U}\{d(\widetilde{m}_i,s)\} \\
&\le& 2\rho \sum_{s\in\S}p(s)\min_{i\in\U}\{d(m^*_i,s)\} \\
&\le& 2\rho J(\mu_*,x(0)).
\end{eqnarray*}
\end{proof}

\vspace{1mm}

In other words, an approximation algorithm which produces factor $\rho$ suboptimal solutions to the $k$-median problem leads to a method for computing factor $2\rho$ suboptimal decentralized policies for the stochastic $k$-server problem. One particularly attractive approximation algorithm for the $k$-median problem is the local search heuristic of \cite{Arya04}. This algorithm is particularly simple to implement and capable of achieving an approximation ratio of $3+\epsilon$ for any $\epsilon > 0$, where there is a tradeoff between computational requirements and approximation ratio.

%%%%%%%%%%%%%%%%%%%%%%%%%%%%%%%%%%%%%%%%%%%%%%%%%%%%%%%%%%%%%%%%%%%%%%%%%%%%%%

\section{Extensions}

In this section we will discuss several extensions of the basic stochastic $k$-server problem and show that results analogous to Theorem \ref{thm:factor_two} can be established.

\subsection{Server-dependent processing costs}

The first extension we consider generalizes the $k$-server model to the case where the servers are not equal in their processing capabilities. In particular, we model the cost of serving a job at location $\widehat{x}$ by server $u$ at location $x_u$ as
\[
r(x,u) = d_u(x_u,\widehat{x})+c_u(\widehat{x}).
\]
The amount of resources consumed (time, fuel, etc.) by moving from location $x_u$ to location $\widehat{x}$ depends on the server, and is modeled by the metric $d_u$ if server $u$ is chosen. Once the server arrives at the service location, an additional cost of $c_u(\widehat{x}) \ge 0$ is incurred when processing the job at location $\widehat{x}$ by server $u$.

As before, decentralized policies partition the state space and assign exactly one server to each partition. We have the following theorem regarding decentralized policies for the case of server-dependent processing costs.

\begin{thm}
\label{thm:server_dep_costs}
For the problem with server-dependent processing costs, there exists a 
decentralized policy $\mu_d$ such that
\[
J(\mu_d,x(0)) \le 2J(\mu_*,x(0))
\]
for all $x(0) \in \X$.
\end{thm}

\vspace{1mm}

\begin{proof}
Similar to the proof of Theorem \ref{thm:factor_two}, we will find a lower bound on $J(\mu_*)$ using Lemma \ref{lem:markov_bounds} with
\[
h_L(x) = \min_{i\in\U}\{d_i(x_i,\widehat{x})+c_i(\widehat{x})\}.
\]
For this choice of $h_L$, we obtain the lower bound
\begin{eqnarray}
\label{eqn:k_median_cost_2}
J(\mu_*,x(0))\bk\ge\bk\min_{m\in\S^k}\left\{\sum_{s\in\S}p(s)\min_{i\in\U}\{d_i(m_i,s)\bk +\bk c_i(s)\} \right\}
\end{eqnarray}
for all $x(0) \in \X$. Note that, unlike the proof of Theorem \ref{thm:factor_two}, the order in which $m_1,\ldots,m_k$ are indexed effects the lower bound in (\ref{eqn:k_median_cost_2}).

Let $m^*$ denote the minimizing $m$ in (\ref{eqn:k_median_cost_2}). The decentralized policy $\mu_d$ divides the set $\S$ into disjoint partitions $\S_1,\ldots,S_k$ where
\[
\S_i = \{s\in\S ~|~ d_i(m^*_i,s)+c_i(s) \le d_j(m^*_j,s)+c_j(s) ~\forall j\}.
\]
We will find an upper bound on $\mu_d$ using Lemma \ref{lem:markov_bounds} with
\[
h_U(x) = 2\min_{i\in\U}\{d_i(m^*_i,\widehat{x})+c_i(\widehat{x})\}+\sum_{i=1}^kd_i(x_i,m^*_i).
\]
Denoting $u_d = \mu_d(x)$, we have
\begin{eqnarray*}
\lefteqn{r(x,u_d) + \Delta_U(x) =} \\ && 2\sum_{s\in\S}p(s)\min_{i\in\U}\{d_i(m^*_i,s)+c_i(s)\} - c_{u_d}(\widehat{x}) + \\ && d_{u_d}(x_{u_d},\widehat{x}) - d_{u_d}(\widehat{x},m^*_{u_d}) - d_{u_d}(x_{u_d},m^*_{u_d}).
\end{eqnarray*}
Since $d_{u_d}$ is a metric,
\[
d_{u_d}(x_{u_d},\widehat{x}) \le d_{u_d}(x_{u_d},m^*_{u_d}) + d_{u_d}(m^*_{u_d},\widehat{x}).
\]
Since $c_{u_d}(\widehat{x}) \ge 0$, we have
\begin{eqnarray*}
J(\mu_d,x(0)) &\le& 2\sum_{s\in\S}p(s)\min_{i\in\U}\{d_i(m^*_i,s)+c_i(s)\} \\
&\le& 2J(\mu_*,x(0)).
\end{eqnarray*}
\end{proof}

\subsection{Multiple requests per period}

Next we consider the case when some fixed number $n \le k$ of requests is generated and must be served in each time step. Specifically, at time step $t$, service is requested at some set of points $\widehat{x}_1(t),\ldots,\widehat{x}_n(t)\in\S$, and exactly $n$ servers must be chosen to service these requests. Here the state at time $t$ is given by $x(t) = (x_1(t),\dots, x_k(t), \widehat{x}_1(t),\ldots,\widehat{x}_n(t))$. Let $u_j(t)$ denote the index of the server chosen to service request $j$. For this case the action at time $t$ is $u(t) = (u_1(t),\ldots,u_n(t))$ and the action space is
\[
\U = \{u\in\{1,\ldots,k\}^n ~|~ u_i\ne u_j \text{ for } i\ne j\}.
\]
At time $t$, a cost of $\sum_{j=1}^nd(x_{u_j(t)}(t),\widehat{x}_j(t))$ is incurred. Server $u_j(t)$ is then relocated to the point $\widehat{x}_j(t)$, and the next set of requests is drawn according to some probability mass function $p:\S^n\rightarrow[0,1]$.

Decentralized policies for this case are a natural extension of the partition policies for the single request case. We will analyze the performance of the decentralized policy $\mu_d$ which is constructed as follows.

\begin{enumerate}
\item Find the $m^*_1,\ldots,m^*_k$ minimizing
\[
\sum_{s\in\S^n}p(s)\min_{u\in\U}\left\{\sum_{j=1}^nd(m_{u_j},s_j)\right\}.
\]
\item Let 
\[
\mu_d(x) = \argmin{u\in\U}\left\{\sum_{j=1}^nd(m^*_{u_j},\widehat{x}_j)\right\}.
\]
\end{enumerate}

\noindent In this policy, the server at point $x_i$ is always associated with the median at point $m_i^*$. When a new batch of requests arrives, each request is matched to one of the medians. No two requests are matched to the same median. If the request at point $\widehat{x}_j$ is matched to the median at point $m_i^*$, then this request is served by the server at point $x_i$. Note that, unlike the single request case, servers may move between partitions associated with several medians. This is because multiple requests may appear in the same partition, and must be served by multiple servers.

Analysis of this case is much like that of the single request case, and is presented in the following theorem.

\begin{thm}
\label{thm:multiple_requests}
The cost of the decentralized policy $\mu_d$ satisfies
\[
J(\mu_d,x(0)) \le 2J(\mu_*,x(0))
\]
for all $x(0) \in \X$.
\end{thm}

\vspace{1mm}

\begin{proof}
The lower bound on $J(\mu_*,x(0))$ is is determined using Lemma \ref{lem:markov_bounds} with
\[
h_L(x) = \min_{u\in\U}\left\{\sum_{j=1}^nd(x_{u_j},\widehat{x}_j)\right\}.
\]
For this choice of $h_L$, we obtain
\begin{eqnarray*}
\label{eqn:multi_k_median_cost}
J(\mu_*,x(0)) &\ge& \min_{m\in\S^k} \left\{ \sum_{s\in\S^n} p(s)\min_{u\in\U}\left\{\sum_{j=1}^nd(m_{u_j},s_j)\right\} \right\} \\
&=&  \sum_{s\in\S^n}p(s)\left(\sum_{j=1}^nd(m^*_{\mu_d(s)_j},s_j)\right)
\end{eqnarray*}
for all $x(0) \in \X$.

The upper bound on $J(\mu_d,x(0))$ is determined using Lemma \ref{lem:markov_bounds} with
\[
h_U(x) = 2\min_{u\in\U}\left\{\sum_{j=1}^nd(m^*_{u_j},\widehat{x}_j)\right\}+\sum_{i=1}^kd(x_i,m^*_i).
\]
Let
\[
x_i(t+1) = \left\{
\begin{array}{cl}
\widehat{x}_j & \text{if } i = \mu_d(x)_j \\
x_i & \text{otherwise}
\end{array}
\right.
\]
For this choice of $h_U$, 
\begin{eqnarray*}
r(x,\mu_d(x))+\Delta_U(x,\mu_d(x)) \bk\bk &=& 2\sum_{s\in\S^n}p(s)\left(\sum_{j=1}^nd(m^*_{\mu_d(s)_j},s_j)\right) \\
&& + \sum_{j=1}^n d(x_{\mu_d(x)_j},\widehat{x}_j) + \sum_{i=1}^kd(x_i(t+1),m^*_i) \\
&& - \sum_{i=1}^kd(x_i,m^*_i) - 2\left(\sum_{j=1}^nd(m^*_{\mu_d(x)_j},\widehat{x}_j)\right) \\
&=& 2\sum_{s\in\S^n}p(s)\left(\sum_{j=1}^nd(m^*_{\mu_d(s)_j},s_j)\right) \\ 
&& + \sum_{j=1}^n\bigl(d(x_{\mu_d(x)_j},\widehat{x}_j) + \sum_{j=1}^n d(m^*_{\mu_d(x)_j},\widehat{x}_j) \\
&& - \sum_{j=1}^nd(x_{\mu_d(x)_j},m^*_{\mu_d(x)_j}) - 2\left(\sum_{j=1}^nd(m^*_{\mu_d(x)_j},\widehat{x}_j)\right) \\
&=& 2\sum_{s\in\S^n}p(s)\left(\sum_{j=1}^nd(m^*_{\mu_d(s)_j},s_j)\right) \\ 
&& \bk+\bk \sum_{j=1}^n \biggl(\bk d(x_{\mu_d(x)_j},\widehat{x}_j) \bk-\bk d(x_{\mu_d(x)_j},m^*_{\mu_d(x)_j}) \bk-\bk d(m^*_{\mu_d(x)_j},\widehat{x}_j)\bk\biggr).
\end{eqnarray*}
Since $d$ is a metric,
\[
d(x_{\mu_d(x)_j},\widehat{x}_j) \le d(x_{\mu_d(x)_j},m^*_{\mu_d(x)_j}) + d(m^*_{\mu_d(x)_j},\widehat{x}_j)
\]
for all $j$. 
Therefore,
\begin{eqnarray*}
J(\mu_d,x(0)) &\le& \inf_{x\in\mathcal{X}}\{r(x,\mu_d(x)) + \Delta_U(x)\} \\
&\le& 2\sum_{s\in\S^n}p(s)\min_{u\in\U}\left\{\sum_{j=1}^nd(m^*_{u_j},s_j)\right\} \\
&\le& 2J(\mu_*,x(0)).
\end{eqnarray*}
\end{proof}

It is worth noting that for the two extensions presented in this section, computing decentralized policies requires solving generalizations of the $k$-median problem. Whether any of the existing approximation algorithms for the $k$-median problem can be extended to these generalizations is not clear, and is a topic for further research.

%%%%%%%%%%%%%%%%%%%%%%%%%%%%%%%%%%%%%%%%%%%%%%%%%%%%%%%%%%%%%%%%%%%%%%%%%%%%%%

\section{Conclusion}

In this paper we presented the \emph{stochastic $k$-server problem}, and showed that a simple decentralized state-feedback policy achieves an average cost within a factor of two of the cost achieved by an optimal state-feedback policy. These results were then extended to several variations of the basic stochastic $k$-server problem.

In this paper, we presented a formulation where the set of possible locations to be served is finite. We have focused on this formulation because low complexity algorithms for computing decentralized policies exist in this case. In fact, it is straightforward to use the results of \cite{Cogill06} to show that the results of this paper hold in infinite bounded metric spaces as well.

%%%%%%%%%%%%%%%%%%%%%%%%%%%%%%%%%%%%%%%%%%%%%%%%%%%%%%%%%%%%%%%%%%%%%%%%%%%%%%

\nocite{Bertsimas93}

\bibliography{stoch_k_server}

\end{document}